\documentclass[11pt,a4paper]{article}
\usepackage[OT1]{fontenc}
\usepackage[latin1]{inputenc}
\usepackage[english]{babel}
\usepackage{amsfonts}
\usepackage{amsthm}
\newtheorem{teo}{Theorem}[section]
\newtheorem{prop}[teo]{Proposition}
\newtheorem{lem}[teo]{Lemma}
\newtheorem{coro}[teo]{Corollary}
\newtheorem{defi}[teo]{Definition}
\theoremstyle{definition}
\newtheorem{rem}[teo]{Remark}

\def\p{{\sf \,p\,}}
\def\h{{\cal H}}

\def\k{{\cal K}}
\def\b{{\cal B}}
\def\a{{\cal A}}

\def\u2{U_2({\cal H})}

\def\u{ {\cal U} }

\begin{document}

\title{\vspace*{0cm} SMOOTH PATHS OF CONDITIONAL EXPECTATIONS}
\date{}
\author{Esteban Andruchow and Gabriel Larotonda}

\maketitle

\abstract{\footnotesize{\noindent  Let $\a$ be a von Neumann algebra with a finite trace $\tau$, represented in $\h=L^2(\a,\tau)$, and let $\b_t\subset \a$ be sub-algebras, for $t$ in an interval $I$ ($0\in I$). Let $E_t:\a\to \b_t$ be the unique $\tau$-preserving conditional expectation. We say that the path $t\mapsto E_t$  is smooth if for every $a\in\a$ and $\xi\in\h$, the map 
$$
I\ni t\mapsto E_t(a)\xi\in \h
$$
is continuously differentiable. This condition implies the existence of the derivative operator
$$
dE_t(a):\h\to\h , \  dE_t(a)\xi=\frac{d}{dt}E_t(a)\xi.
$$
If this operator verifies the additional boundedness  condition,
$$
\int_J \|dE_t(a)\|_2^2 d t\le C_J\|a\|_2^2,
$$
for any closed bounded sub-interval $J\subset I$, and $C_J>0$ a constant depending only on $J$, then the algebras $\b_t$ are $*$-isomorphic. More precisely, there exists a curve $G_t:\a\to\a$, $t\in I$ of  unital, $*$-preserving linear isomorphisms which intertwine the expectations,
$$
G_t\circ E_0=E_t\circ G_t.
$$ The curve $G_t$ is weakly continuously differentiable. Moreover, the intertwining property in particular implies that $G_t$ maps $\b_0$ onto $\b_t$. We show that this restriction is  a multiplicative  isomorphism.
}}

\setlength{\parindent}{0cm} 

\section{Introduction}
 Let $\a$ be a von Neumann algebra with a finite faithful and normal trace $\tau$, and suppose $\a$ acting on its standard Hilbert space $\h=L^2(\a,\tau)$. We shall assume that for each $t\in I$ ($0\in I$), there is a von Neumann sub-algebra $\b_t\subset \a$, and we shall denote by $E_t:\a\to\b_t$ the unique $\tau$-invariant conditional expectation. We  regard $t\mapsto E_t$ as a curve, and require smoothness in the following sense: for each $a\in\a$ and $\xi\in\h$,
$I\ni t\mapsto E_t(a)\xi\in\h$ is continuously differentiable. This paper is a sequel to \cite{weaksmooth}, where a similar matter is treated with more strict hypothesis. In \cite{weaksmooth} we considered a stronger smoothness condition, namely, that for each $a\in\a$, the map $I\ni t\mapsto E_t(a)\in\a$ is continuously differentiable (in norm).

The current regularity  assumption  on $E_t$ implies the existence of the bounded derivative operator, for each $t\in I$ and $a\in\a$
$$
dE_t(a):\h\to\h , \ \ dE_t(a)\xi=\frac{d}{d t}E_t(a)\xi.
$$
Therefore a curve of possibly unbounded symmetric operators $dE_t$ is defined in $\h$, with common domain $\a\subset\h$.
We shall make the following assumption on $dE$:
\begin{equation}\label{hipotesis}
\int_J\|dE_t(a)\|_2^2 dt\le C_J\|a\|_2^2
\end{equation}
for all $a\in\a$, and every closed bounded interval $J\subset I$ (the constant depends only on $J$).

With these assumptions, we prove that there exists a curve $I\ni t\mapsto G_t$ of linear isomorphisms $G_t:\a\mapsto \a$ with the following properties:
\begin{enumerate}
\item
 For each $a\in \a$, the curve $I\ni t \to G_t(a)\in \a\subset \h$ is weakly continuously differentiable, with $G_0=Id$.
  \item
 The maps $G_t$ are unital and $*$-preserving.
 \item
 For each $t\in J_0$, 
 $$
 G_t E_0 G_t^{-1}=E_t.
 $$
 \item
 The last formula implies that $G_t$ maps $\b_0$ onto $\b_t$. The restriction
 $$
 G_t|_{\b_0}:\b_0\to \b_t
 $$
 is a $*$-isomorphism.
\item
The linear isomorphisms $G_t:\a\to\a$ are $\|\ \|_2$-isometric, therefore they extend to unitary operators $U_t$ acting in $\h$, which preserve $\a$ ($U_t(\a)=\a$).
 \end{enumerate}

A similar result was obtained in \cite{weaksmooth} with the already noted stronger assumption. In both contexts, the maps $G_t$ appear as propagators of the linear differential equation
\begin{equation}\label{trasporte}
\left\{ \begin{array}{l} \dot{\alpha}(t)=dE_t(E_t(\alpha(t)))-E_t(dE_t(\alpha(t))) \\ \alpha(s)=a, \end{array} \right. 
\end{equation}
for $a\in\a$.
In the present context, our hypothesis does not guarantee that the linear operators $[dE,E]$ of this equation are bounded, nor that they vary continuously. Therefore our first task is to show that with the current assumptions (particularly \ref{hipotesis}), this equation has existence and uniqueness of {\it weak} solutions. This is done in section 3. In section 2 we state the basic properties of the operator $dE$. In section 4 we prove the existence and properties of the maps $G_t$. In section 5 we consider the example when the expectations $E_t$ are given by a curve of systems of projections $p_1(t),p_2(t),...$ in $\a$ (i.e. curves of pairwise orthogonal projections which sum up to $1$), and examine when our hypothesis are verified.

\section{Curves of expectations}
As we stated above, we shall consider $\a$ represented in the standard space $\h=L^2(\a,\tau)$, and also regard elements of $a$ as elements in $\h$. We shall denote by $\| \ \|_\infty$ the norm of $\a$, and by $\|\ \|_2$ the norm of $\h$. 

\begin{lem}
For each $a\in\a$ and $t\in I$, the linear operator $dE_t(a)$ defined in the previous section is bounded, its adjoint is $dE_t(a^*)$.
\end{lem}
\begin{proof}
Note that both $dE_t(a)$ and $dE_t(a^*)$ are defined in the whole space $\h$  by hypothesis. If $x,y\in \a$, regarded as a dense subspace of  $\h$,
$$
<dE_t(a)x,y>=\frac{d}{d t}<E_t(a)x,y>=\frac{d}{d t}\tau(y^*E_t(a)x)=\frac{d}{d t}\tau((E_t(a^*)y)^*x)
$$
$$
=\frac{d}{d t}<x,E_t(a^*)y>=<x,dE_t(a^*)y>.
$$
By the closed graph theorem, it follows that $dE_t(a)$ is bounded, and that $dE_t(a^*)$ is its adjoint.
\end{proof}
Next let us show that the derivative of $E_t$ defines also a map on $\a$. 
\begin{lem}
Let $a\in \a$, then for each $t\in I$, $dE_t(a)\in \a$. 
\end{lem}
\begin{proof}
Let $T\in\b(\h)$ belong to the commutant of $\a$. If $\xi,\eta\in\h$,
$$
<dE_t(a)T\xi,\eta>=\frac{d}{dt}<E_t(a)T\xi,\eta>=\frac{d}{d t}<TE_t(a)\xi,\eta>=\frac{d}{d t}<E_t(a)\xi,T^*\eta>
$$
$$
=<dE_t(a)\xi,T^*\eta>=<TdE_t(a)\xi,\eta>,
$$
i.e. $dE_t(a)\in\a$.
\end{proof}
The correspondence $dE_t:\a\to\a$ is apparently linear, and $*$-preserving. Let us verify that it is bounded as an operator acting in $(\a,\| \ \|_\infty)$. 
\begin{prop}\label{acotacion dE}
For each $t\in I$, the map $dE_t:(\a,\|\ \|_\infty)\to (\a,\|\ \|_\infty)$, $a\mapsto dE_t(a)$, is linear, $*$-preserving and bounded. Moreover, for any closed and bounded sub-interval $J\subset I$, the norms of the operators $dE_t:(\a,\|\ \|_\infty)\to (\a,\|\ \|_\infty)$, denoted $\|dE_t\|_{\infty,\infty}$, are uniformly bounded for $t\in J$.
\end{prop}
\begin{proof}
Let us prove that the graph of $dE_t$ is closed. Let $a_n,a, b\in\a$ such that $\|a_n-a\|_\infty \to 0$ and $\|dE_t(a_n)-b\|_\infty\to 0$. First note that if $x,y\in \a$, then
$$
\tau(dE_t(x)y)=\tau(xdE_t(y)).
$$
Indeed, by the invariance of $E_t$ and $\tau$, 
$$
\tau(E_t(x)y)=\tau(E_t(E_t(x)y))=\tau(E_t(x)E_t(y))=\tau(E_t(xE_t(y)))=\tau(xE_t(y)).
$$
Then 
$$
\tau(dE_t(x)y)=<dE_t(x),y^*>=\frac{d}{d t}<E_t(x),y^*>=\frac{d}{d t}\tau(E_t(x)y)=\frac{d}{d t}\tau(xE_t(y)),
$$
which by the same argument equals $\tau(xdE_t(y))$.
Therefore, for any $x\in \a$,
$$
\tau(bx)=\lim_{n\to\infty}\tau(dE_t(a_n)x)=\lim_{n\to\infty}\tau(a_ndE_t(x))=\tau(a\ dE_t(x))=\tau(dE_t(a)x).
$$
It follows that $dE_t(a)=b$, and therefore $dE_t$ is bounded. 

Consider now a closed bounded sub-interval $J\subset I$. Fix $a\in \a$ and $\xi\in\h$. Since by hypothesis the map $t\mapsto E_t(a)\xi$ is continuously differentiable, it follows that there exists a constant $C_{J,a,\xi}$
 such that
 $$
 \|dE_t(a)\xi\|_2\le C_{J,a,\xi} \ \ \hbox{ for all } t\in J.
 $$
 By the uniform boundedness principle in the Banach space $(\h, \| \ \|_2)$, it follows that there exists a constant $C_{J,a}$ such that
 $$
 \|dE_t(a)\|_\infty\le C_{J,a} \ \ \hbox{ for all } t\in J.
 $$
 Again by the uniform boundedness principle, this time in the Banach space $(\a,\|\ \|_\infty)$, it follows that there exists a constant $C_J$ such that
 $$
 \|dE_t\|_{\infty,\infty}\le C_J \ \ \hbox{ for all } t\in J.
 $$
 \end{proof} 
 We emphasize that $dE_t$ may be an unbounded operator in $\h$, with domain $\a$.
 \begin{rem}
 The assumption that $I\ni t\mapsto E_t(a)\xi\in \h$ is continuously differentiable implies  that $t\mapsto E_t(a)\in\h$ is continuously differentiable. Indeed, it suffices to take $\xi=1\in\a$.

 \end{rem}
 We shall need the following elementary fact.
 \begin{lem}\label{derivada del producto}
 For $h\in[-\delta,\delta]$, let $b_h,b\in\a$ such that $\|b_h-b\|_2\to 0$ as $h\to 0$. Then 
 $$
 \|E_{t+h}(b_h)- E_t(b)\|_2\to 0 \ \ \hbox{ as } h\to 0.
 $$
 \end{lem}
 \begin{proof}
 Note  that
 $$
 \|E_{t+h}(b_h)-E_t(b)\|_2\le \|E_{t+h}(b_h)-E_{t+h}(b)\|_2+\|E_{t+h}(b)-E_t(b)\|_2.
 $$
 The second term clearly tends to $0$. Since the expectations $E_t$ are $\tau$-invariant, they are contractive for the $\|\ \|_2$-norm. Therefore the first term is bounded by $\|b_h-b\|_2$.
\end{proof}
We shall use the following formula thoroughly.
\begin{prop}\label{codiagonal}
For any $a\in\a$ and any $t\in I$,
$$
dE_t(E_t(a))+E_t(dE_t(a))=dE_t(a).
$$
\end{prop}
\begin{proof}
Note that
$$
\frac1h \{E_{t+h}(a)-E_t(a)\}=\frac1h \{E_{t+h}(E_{t+h}(a))-E_t(E_t(a))\}
$$
$$
=E_{t+h}(\frac1h\{E_{t+h}(a)-E_t(a)\})+\frac1h\{E_{t+h}(E_t(a))-E_t(E_t(a))\}.
$$
The second term tends to $dE_t(E_t(a))$ in the $2$-norm.  The first term tends to $E_t(dE_t(a))$ in the $2$-norm by the above Lemma, which proves the formula.
\end{proof}
\section{The transport equation}
Under the current assumptions we shall examine existence and uniqueness of solutions of the linear differential equation below, which we shall call the transport equation (\ref{trasporte})
$$
\left\{ \begin{array}{l} \dot{\alpha}(t)=dE_t(E_t(\alpha(t)))-E_t(dE_t(\alpha(t))) \\ \alpha(s)=a, \end{array} \right. 
$$
where $a\in \a$. We shall be looking for solutions $\alpha(t)$ with values in $\a$, which are differentiable as $\h$-valued maps in the weak sense. That is, $t\mapsto <\alpha(t),\xi>$ is differentiable, and its derivative verifies
$$
\frac{d}{d t}<\alpha(t),\xi>=<dE_t(E_t(\alpha(t)))-E_t(dE_t(\alpha(t))),\xi>,
$$
for all $\xi\in\h$.

Note that the classical results on linear differential equations in Banach spaces (for instance, \cite{skrein,reedsimon2}) do not apply. The linear operators  $[dE_t,E_t]$ need not be continuous in the parameter $t$ as operators in the Banach space $\a$, nor they need to be bounded as operators in $\h$ (with common domain $\a$), or even closed operators. This seems to be a mixed terrain, where both considerations with the non equivalent norms $\|\  \|_\infty$ and $\| \ \|_2$ play a role. We shall show existence and uniqueness of solutions mimicking carefully Picard's method of succesive approximations, under the assumption of the following Hypothesis (\ref{hipotesis}):
$$
\int_J\|dE_t(a)\|_2^2 dt\le C_J\|a\|_2^2
$$
for all $a\in\a$, and every closed bounded interval $J\subset I$ (the constant depends only on $J$).

Note that this hypothesis trivially holds if $dE$ is bounded in the $2$-norm $\| \ \|_2$. Indeed, this holds by the uniform boundedness principle.

We shall mainly be involved with the properties of the operator $H_t=[dE_t,E_t]$. Note that $H_t(\a)\subset \a$. Also it is clear that $H_t$ is anti-symmetric in $\a$: if $x,y\in \a$ then 
$$
<H_t(x),y>=<dE_t(E_t(x)),y>-<E_t(dE_t(x)),y>
$$
$$
=<x,E_t(dE_t(y)>-<x,dE_t(E_t(y))>=-<x,H_t(y)>.
$$
Also it is apparent that for each fixed $x\in\a$, $t\mapsto H_t(x)\in \h$ is continuous.

The following result will be needed. It is not supposed in the next Lemma that $E_t$ verifies Hypothesis (\ref{hipotesis}).
\begin{lem}
Let $f:I\to \a$ be uniformly $\|\ \|_\infty$-bounded on closed bounded sub-intervals of $I$, and weakly continuous when regarded as an $\h$-valued map, i.e.
\begin{enumerate}
\item 
For every closed bounded $J\subset I$ there exists a constant $C_J$ such that $\|f(t)\|_\infty\le C_J$ for all $t\in J$.
\item
For every $\xi\in\h$, the map $t\mapsto <f(t),\xi>$ is continuous.
\end{enumerate}
Then the map $t\mapsto H_t(f(t))$ takes values in $\a$, is  weakly continuous as an $\h$-valued map, and is uniformly $\|\ \|_\infty$-bounded on closed bounded intervals as an $\a$-valued map.
\end{lem}
\begin{proof}
First pick $x\in\a$. Then $g_x(t)=<H_t(f(t)),x>=-<f(t),H_t(x)>$.
Thus 
$$
g_x(t+h)-g_x(t)=-<f(t+h),H_{t+h}(x)>+<f(t),H_t(x)>
$$
$$
=<f(t+h),H_t(x)-H_{t+h}(x)>+<f(t+h)-f(t),H_t(x)>.
$$
The second term tends to $0$ as $h\to 0$. By the Cauchy-Scwarz inequality, the first term is bounded by
$$
\|f(t+h)\|_2\|H_{t+h}(x)-H_t(x)\|_2.
$$
This expression also tends to $0$, as $h\to 0$, because $f$ is $\| \ \|_\infty$ bounded (and therefore also $\| \ \|_2$ bounded). Let $\xi\in\h$ and pick $x\in \a$ such that $\|\xi-x\|_2<\epsilon$. Then if $g_\xi(t)=<H_t(f(t)),\xi>$,
$$
g_\xi(t+h)-g_\xi(t)=<H_{t+h}(f(t)),\xi-x>+g_x(t+h)-g_x(t)+<H_t(f(t)),x-\xi>.
$$
If $h\to 0$, the middle term tends to $0$. Again, by the  Cauchy-Scwarz inequality, the first term is bounded by
$$
\|H_{t+h}(f(t+h))\|_2\|\xi-x\|_2\le \|H_{t+h}(f(t+h))\|_\infty \|\xi-x\|_2\le \epsilon\|H_{t+h}\|_{\infty,\infty}\|f(t+h)\|_\infty.
$$
For small $h$ (e.g. $|h|\le \delta$ such that $J=[t-\delta,t+\delta]\subset I$), both factors above are uniformly bounded. For instance $\|H_t\|_{\infty,\infty}\le 2\|dE_t\|_{\infty,\infty}$, and then use Proposition \ref{acotacion dE}.
The third term is dealt similarly. This proves the weak continuity of $t\mapsto H_t(f(t))\subset \h$.

Local boundedness in $\|\ \|_\infty$ is straightforward: $\|H_t(f(t))\|_\infty\le 2\|dE_t\|_{\infty,\infty} \|f(t)\|_\infty$.
\end{proof}

Fix $a\in\a$ and $s$ in the interior of $I$. For each $t\in I$, consider the following sequence of functions $S^{a,s}_n(t)=S_n(t)$:
\begin{defi}\label{sucesion de los sn}
$$
S_0(t)=a \ , \ \ S_1(t)=a+{\bf weak}\int_s^t H_u(a)du , \ \hbox{ and } \ \  S_{n+1}(t)=a+{\bf weak}\int_s^t H_u(S_n(u)) du,
$$
where ${\bf weak}\int$ stands for the {\it weak integral}, i.e. for each $\xi\in \h$, ${\bf weak}\int_J f(u) du$ is given by
$$
<{\bf weak}\int_J f(u) du,\xi>=\int_J<f(u),\xi>du.
$$
\end{defi}
First we must show that $S_n(t)$ is well defined.
\begin{prop}
For any fixed $a\in\a$ and $s$ in the interior of $I$, the maps $S_n(t)$, $t\in I$ are well defined. They take values in $\a$. Regarded as $\a$-valued functions, they are uniformly bounded on  closed bounded sub-intervals of $I$. Regarded as $\h$-valued functions, they are weakly continuous.
\end{prop}
\begin{proof}
This is proved by induction. Clearly $S_0$ takes values in $\a$, is $\| \ \|_\infty$-bounded uniformly bounded on closed bounded intervals,  and is $\h$-weakly continuous. Suppose that $S_n$ verifies these conditions. By  the above lemma, the map $t\mapsto H_t(S_n(t))$ is $\h$-weakly continuous and $\| \ \|_\infty$-bounded. Therefore, it only remains to be verified that it takes values in $\a$. The weak integral $\int_s^t H_u(S_n(u)) du$ is the weak limit of its Riemann sums $\sum_jH_{u_j}(S_n(u_j))(u_j-u_{j-1})$, which are linear combinations of elements of $\a$, and thus lie in $\a$. Moreover 
$$
\|\sum_jH_{u_j}(S_n(u_j))(u_j-u_{j-1})\|_\infty\le \sum_j\|H_{u_j}(S_n(u_j))\|_\infty(u_j-u_{j-1}).
$$
Each term $\|H_{u_j}(S_n(u_j))\|_\infty$ is uniformly bounded in the interval $[s,t]$. Therefore the Riemann sums are uniformly $\|\ \|_\infty$-bounded. Therefore the weak limit of these sums lies in $\a$.
\end{proof}
For the next result we need Hypothesis (\ref{hipotesis})
\begin{prop}\label{acotacion}
Fix $s_0 \le t_0$ in $I$ and $a\in \a$, and consider $S_n(t)=S_n^{s_0,a}(t)$. Assume that Hypothesis (\ref{hipotesis}) holds: $\int_{s_0}^{t_0}\|dE_s(b)\|_2^2ds\le C\|b\|_2^2$ (where $C=C_{[s_0,t_0]}$).
Then for all $t\in[s_0,t_0]$,
$$
\|S_{n+1}(t)-S_n(t)\|_2\le C^{1/2}\sqrt{t-s_0} \sup_{u\in[s_0,t]} \|S_n(u)-S_{n-1}(u)\|_2.
$$
\end{prop}
\begin{proof}
Pick $b\in\a$. Then
$$
|<S_{n+1}(t)-S_n(t),b>|=|\int_{s_0}^t<H_u(S_n(u))-H_u(S_{n-1}(u)), b>du|
$$
$$
=|\int_{s_0}^t<S_n(u)-S_{n-1}(u), H_u(b)>du|\le \int_{s_0}^t |<S_n(u)-S_{n-1}(u),H_u(b)>| du
$$
$$
\le \sup_{u\in [s_0,t]}\|S_n(u)-S_{n-1}(u)\|_2\int_{s_0}^t \|H_u(b)\|_2 du.
$$
By H\"older's inequality
$$
\int_{s_0}^t \|H_u(b)\|_2 du\le \{\int_{s_0}^t \|H_u(b)\|_2^2 du\}^{1/2}\sqrt{t-s_0}.
$$
Recall that $H_u(b)=dE_u(E_u(b))-E_u(dE_u(b))$. Using the formula  in Proposition \ref{codiagonal}, $dE_u(b)=dE_u(E_u(b))+E_u(dE_u(b))$, one obtains that
$$
H_u(b)=dE_u(b)-2E_u(dE_u(b))=(1-2E_u)(dE_u(b)).
$$
Note that $E_u$ is (or rather, extends to) a self
adjoint projection in $\h$. Therefore $1-2E_u$ is a symmetry, i.e. a selfadjoint unitary operator. In particular, it is $\|\ \|_2$-isometric. Therefore
$$
\|H_u(b)\|_2=\|(1-2E_u)(dE_u(b))\|_2=\|dE_u(b)\|_2.
$$
Then (using Hypothesis (\ref{hipotesis}))
$$
|<S_{n+1}(t)-S_n(t),b>|\le \sup_{u\in [s_0,t]}\|S_n(u)-S_{n-1}(u)\|_2\{\int_{s_0}^t \|dE_u(b)\|_2^2du \}^{1/2}\sqrt{t-s_0}
$$
$$
\le \sup_{u\in [s_0,t]}\|S_n(u)-S_{n-1}(u)\|_2C^{1/2}\|b\|_2\sqrt{t-s_0}.
$$
Taking supremum over $b\in\a$ with $\|b\|_2=1$ proves the inequality.
\end{proof}
\begin{coro}
Fix $s_0\in I$ and $a\in \a$.  If Hypothesis (\ref{hipotesis}) holds, then there exists $t_0\in I$, $s_0<t_0$, such that the sequence $S_n^{s_0,a}(t)=S_n(t)$ converges uniformly in the norm $\| \ \|_2$, in the interval $[s_0,t_0]$,  to a function $S(t)$. This function $S(t)$ takes values in $\a$, is uniformly $\|\ \|_\infty$-bounded, and  weakly continuously differentiable as an $\h$-valued map. Moreover, for $t\in [s_0,t_0]$ and $\xi\in\h$,
$$
<S(t),\xi>=<a,\xi>+\int_{s_0}^t <H_s(S(s)),\xi> d s.
$$
\end{coro}
\begin{proof}
Pick $t_0$ such that $k_0=C^{1/2}\sqrt{t_0-s_0}<1$, where $C$ is the constant in the above Proposition. Then, if $t\in[s_0,t_0]$,
$$
\|S_{n+1}(t)-S_n(t)\|_2\le C^{1/2}\sqrt{t-s_0}\sup_{u\in[s_0,t]}\|S_n(u)-S_{n-1}(u)\|_2
$$
$$
\le C^{1/2}\sqrt{t_0-s_0}\sup_{u\in[s_0,t_0]}\|S_n(u)-S_{n-1}(u)\|_2=k_0 \sup_{u\in[s_0,t_0]}\|S_n(u)-S_{n-1}(u)\|_2.
$$
Then 
$$
\sup_{t\in[s_0,t_0]}\|S_{n+1}(u)-S_n(t)\|_2\le k_0\sup_{t\in[s_0,t_0]}\|S_n(t)-S_{n-1}(t)\|_2.
$$
It follows, by a well-known argument, that $S_n(t)$ converges in $\h$ to a function $S(t)$, uniformly in $[s_0,t_0]$.
The maps $S_n(t)$ are $\a$-valued and uniformly $\|\ \|_\infty$-bounded in $[s_0,t_0]$, therefore $S(t)$ is also $\a$-valued, and uniformly $\|\ \|_\infty$-bounded. Note that it is weakly continuous as an $\h$-valued map: if $\xi\in\h$, then $<S(t+h)-S(t),\xi>$ equals 
$$
<S(t+h)-S_n(t+h),\xi>+<S_n(t+h)-S_n(t),\xi>+<S_n(t)-S(t),\xi>.
$$
and the proof follows by a typical $\epsilon/3$ argument. Finally, by construction, for any $x\in \a$
$$
<S_{n+1}(t),x>=<a,x>+\int_{s_0}^t <H_u(S_n(u)),x> d u=<a,x>-\int_{s_0}^t <S_n(u),H_u(x)> d u.
$$
Note that $<S_n(u),H_u(x)>$ tends uniformly to $<S(u),H_u(x)>$ in the interval $[s_0,t_0]$. Indeed,
$$
|<S_n(u),H_u(x)>-<S(u),H_u(x)>|\le \|S_n(u)-S(u)\|_2\|H_u(x)\|_2
$$
$$
\le \|S_n(u)-S(u)\|_2\|H_u(x)\|_\infty,
$$
where, as seen before, $\|H_u(x)\|_\infty$ is uniformly bounded in $[s_0,t_0]$. Therefore, in the expression above, 
taking limit $n\to \infty$, one obtains
$$
<S(t),x>=<a,x>+\int_{s_0}^t <H_u(S(u)),x>
$$
for all $x\in \a$. By density, it follows that
$$
<S(t),\xi>=<a,\xi>+\int_{s_0}^t <H_u(S(u)),\xi>
$$
for all $\xi\in \h$. In particular, this implies that $S(t)$ is weakly continuously differentiable as an $\h$-valued map.
\end{proof}

The next step is to extend this weak solution. Fix a closed bounded interval $J_0\subset I$, and let $C=C_{J_0}$ be the constant in the inequality of Hypothesis (\ref{hipotesis}) for this sub-interval. If $s_0\in J_0$, then the length of the interval $[s_0,t_0]$ on which a solution is defined depends only on this constant $C$. It does not depend on the initial condition $a$. It follows that one can glue solutions in a standard fashion, to obtain a solution $S(t)$ defined in the whole sub-interval $J_0$. Uniqueness of solutions follows. Indeed, suppose that $S_1,S_2$ are two solutions with $S_1(s)=S_2(s)$.
Then
$$
S_i(t)=a+{\bf weak}\int_s^t H_u(S_i(u)) du \, \ \  i=1,2.
$$
Thus, as in Proposition \ref{acotacion},  
$$
\|S_1(t)-S_2(t)\|_2\le C_{J_0}^{1/2}\sqrt{t-s} \sup_{u\in[s,t]} \|S_1(t)-S_2(t)\|_2.
$$
Then $S_1$ and $S_2$ coincide up to time $t$ such that $|t-s|<1/C_{J_0}$. Note  that this constant does not depend on $s$. It follows that $S_1$ and $S_2$ coincide in $J_0$. Clearly this holds on any closed bounded sub-interval $J_0\subset I$.

Let us summarize these  results.
\begin{teo}\label{existenciayunicidad}
Suppose that Hypothesis (\ref{hipotesis}) holds.
Let $a\in\a$,   Then there exists a map $\alpha_s(t)$, which is $\a$-valued, uniformly $\| \ \|_\infty$-bounded on closed bounded subintervals of $I$, and  weakly continuously differentiable as an $\h$-valued function,  which is the unique (weak) solution of the transport  equation (\ref{trasporte})
$$
\left\{ \begin{array}{l} \dot{\alpha}(t)=[dE_t,E_t](\alpha(t)) \\ \alpha(s)=a. \end{array} \right.
$$
\end{teo}
\begin{rem}\label{posteorema}
For $s,t\in I$, denote by $G_{t,s}$ the propagator of the transport equation, i.e.
$$
G_{t,s}:\a\to \a \ , \ \ G_{t,s}(a)=\alpha_s(t),
$$
where $\alpha_s$ is the solution of (\ref{trasporte}) with $\alpha_s(s)=a$. The propagator has the following properties:
\begin{enumerate}
\item
 $G_{t,s}$ is isometric for the $\|\ \|_2$ norm: $\|G_{t,s}(a)\|_2=\|a\|_2$.
\item
For each $a\in\a$, $G_{t,s}(a)$, as an $\h$-valued map, is weakly continuously differentiable in the parameter $t$, and  continuous in the parameter $s$. 
\item
$G_{s,s}(a)=a$, for all $a\in\a$. 
\item
$G_{t,s}G_{s,r}=G_{t,r}$.
\end{enumerate}
To prove the first assertion, put $\alpha_s(t)=G_{t,s}(a)$, ($\alpha_s(s)=a$), then
 $$ 
 \frac{d}{dt}<G_{t,s}(a),G_{t,s}(a)>=<H_t(\alpha_s(t)),\alpha_s(t)>+<\alpha_s(t),H_t(\alpha_s(t))>=0.
 $$
 Here we use the fact that the product rule holds for weak solutions because they are uniformly $\|\ \|_\infty$-bounded, and also that $H_t=[dE_t,E_t]$ is anti-symmetric. Therefore 
 $$
 \|G_{t,s}(a)\|_2^2=\|G_{s,s}(a)\|_2^2=\|a\|_2^2.
 $$
The third and fourth assertions are apparent. 
To prove the second, use the fourth:
$$
G_{t,s+h}(a)-G_{t,s}(a)=G_{t,s}(G_{s,s+h}(a)-a).
$$
And then, for $b\in \a$,
$$
<G_{t,s+h}(a)-G_{t,s}(a),b>=<G_{s,s+h}(a)-a,G_{t,s}^*(b)>
$$
$$
=\int_s^{s+h} <H_u(G_{u,s+h}(a)-a),G_{t,s}^*(b)> du.
$$
For $|h|<\delta$ such that $[s-\delta,s+\delta]\subset I$ there exists a constant $D$ such that $\|dE_u\|_{\infty,\infty}\le D$. Then
$$
\|H_u(G_{u,s+h}(a)-a)\|_2=\|dE_u(G_{u,s+h}(a)-a)\|_2\le \|dE_u(G_{u,s+h}(a)-a)\|_\infty 
$$
$$
\le D\|G_{u,s+h}(a)-a\|_\infty,
$$ 
which is uniformly bounded for such $h$, by a constant $D'$. Therefore
  $$
 |<G_{t,s+h}(a)-G_{t,s}(a),b>|\le |\int_s^{s+h}|<H_u(G_{u,s+h}(a)-a),G_{t,s}^*(b)>|  \  \ du |
 $$
 $$
  \le  |\int_s^{s+h}\|H_u(G_{u,s+h}(a)-a)\|_2 \|b\|_2 du| \le D' |h| \|b\|_2.
 $$
 Taking supremum over $b\in \a$ with $\|b\|_2=1$, one has
 $$
\|G_{t,s+h}(a)-G_{t,s}(a)\|_2 \le D'|h|.
$$
Note that one obtains more than continuity in the parameter $s$.

In particular, these facts  imply that the map
\begin{equation}\label{el mapa G}
G_t:\a\to \a \ , \ \ G_t:=G_{t,0}
\end{equation}
is invertible, its inverse is $G_t^{-1}=G_{0,t}$.
\end{rem}

\section{The propagators as intertwiners}

In this section we show that the linear isomorphisms
 $G_t$ intertwine the expectations: 
 $$
 G_t\circ E_0\circ G_t^{-1}=E_t.
 $$
 To this effect, the following result is needed.
\begin{prop}\label{solucionproyectada}
Let $\alpha(t)$, $t\in I$ be a (weak) solution of the transport equation (\ref{trasporte}). Then the map $E_t(\alpha(t))$ is also a solution. In particular, if at any given instant $t_0\in I$ one has that $\alpha(t_0)\in\b_{t_0}$, then $\alpha(t)\in \b_t$ for all $t\in I$.
\end{prop}
\begin{proof}
First we must show that $\beta=E(\alpha)$ is $\a$-valued, $\|\ \|_\infty$-bounded and weakly continuously differentiable as an $\h$-valued function. The first fact is apparent. The second: $\|E_t(\alpha(t))\|_\infty\le \|\alpha(t)\|_\infty$. The third: if $\xi \in \h$
$$
\frac1h<\beta(t+h)-\beta(t),\xi>=<E_{t+h}(\frac{\alpha(t+h)-\alpha(t)}{h}),\xi>+<(\frac{E_{t+h}-E_t}{h})(\alpha(t)),\xi>.$$
The second term tends to $<dE_t(\alpha(t)),\xi>$ as $h\to 0$, by definition. For the first term we can aply Lemma \ref{derivada del producto}, and it follows that it tends to $<E_t(\dot{\alpha}(t)),\xi>$. Then $E(\alpha)$ is weakly differentiable, and its derivative is $dE(\alpha)+E(\dot{\alpha})$, which is weakly continuous.
Let us verify that $E(\alpha)$ is a solution:
$$
\frac{d}{dt}E(\alpha)=dE(\alpha)+E(\dot{\alpha})=dE(\alpha)+E(dE(E(\alpha)))-E(E(dE(\alpha))).
$$
Recall  from Lemma \ref{codiagonal} that $dE=dE(E) +E(dE)$, which in particular implies that 
$$
E ( dE ) E=0.
$$
 Then the expression above equals
$$
dE(\alpha)-E(dE(\alpha))=dE(E(\alpha)).
$$
On the other hand
$$
[dE,E](E(\alpha))=dE(E(E\alpha))-E(dE(E(\alpha)))=dE(E(\alpha)).
$$
The last assertion follows by uniqueness of solutions.
\end{proof}
Our main result follows:
\begin{teo}\label{propagador}
Let $E_t:\a\to \b_t\subset \a$, $t\in I$ be a curve of trace invariant conditional expectations, such that for each $x\in\a$ and $\xi\in\h$, the $\h$-valued curve $E_t(x)\xi$ is continuously differentiable. Suppose also that $E_t$ verifies Hypothesis (\ref{hipotesis}), i.e. for each closed bounded subinterval $J\subset I$, 
$$
\int_J \|dE_t(a)\|_2^2 d t \le C_J \|a\|_2^2.
$$
Then  the curve of propagators $G_t:\a\to \a$, $t\in I$, verifies:
\begin{enumerate}
\item
 For each $a\in \a$, the curve $I\ni t \to G_t(a)\in \a\subset \h$ is weakly continuously differentiable, with $G_0=Id$.
  \item
 The maps $G_t$ are unital and $*$-preserving.
 \item
 For each $t\in I$, 
 $$
 G_t E_0 G_t^{-1}=E_t.
 $$
 \end{enumerate}

 \end{teo}
 \begin{proof}
 The first assertion is apparent: $G_t(a)$ is a weak solution of the transport equation. 
 
 Since $E_t(1)=1$ for all $t$, $dE_t(1)=0$, and therefore $H_t(1)=0$. Therefore $\alpha(t)=1$ for all $t$ is a solution, i.e. $G_t(1)=1$. The maps $E_t$ are also $*$-preserving: $E_t(a^*)=E_t(a)^*$, therefore also $dE_t(a^*)=dE_t(a)^*$ and $H_t(a^*)=H_t(a)^*$. Therefore if $\alpha(t)$ is a solution, then also $\alpha^*(t)$ is a solution, and thus $G_t(a^*)=G_t(a)^*$. For the last assertion, note that by the above Proposition, $E_t(G_t(a))$ is a solution. Clearly also $G_t(E_0(a))$ is a solution. At $t=0$, they take the values $E_0(G_0(a))=E_0(a)$ and $G_0(E_0(a))=E_0(a)$, therefore $E_t(G_t(a))=G_t(E_0(a))$ for all $t\in I$. 
 \end{proof}
 \begin{rem} Under the hypothesis of the above theorem, the first assertion in Remark \ref{posteorema} implies that the propagators  $G_t:\a\to \a$ can be extended to unitary operators $U_t$ acting in $\h$. Clearly they preserve $\a\subset \h$: $U_t(\a)\subset \a$. Moreover, if $e_t$ denotes the extension of $E_t$ to an operator in  $\h$, in fact a selfadjoint projection, the last assertion implies that these projections are unitarily equivalent, more precisely
 $$
 U_te_0U_t^*=e_t \ , \ \ t\in I.
 $$
 \end{rem}
 The identity $G_t E_0 G_t^{-1}=E_t$ of the above theorem, in particular implies that $G_t$ maps $\b_0$ onto $\b_t$. Our next result shows that this restriction is a multiplicative $*$-isomorphism.

 \begin{teo}\label{multiplicativo}
 Assume Hypothesis (\ref{hipotesis}). Then for each $t\in I$, the map $\theta_t:=G_t|_{\b_0}:\b_0\to \b_t$ is a multiplicative $*$-isomorphism.
 \end{teo}
 \begin{proof}
 The above indentity clearly implies that $\theta_t(\b_0)=\b_t$. Also it is clear that $\theta_t$ is linear, $*$-preserving and bijective. Thus it only remains to prove that it is multiplicative. Let $a,b\in\b_0$, and denote by $\alpha$ and $\beta$ the solutions of the transport equation with $\alpha(0)=a$ and $\beta(0)=b$. Note that Proposition \ref{solucionproyectada} implies that both $\alpha(t),\beta(t)\in\b_t$, i.e. $E_t(\alpha(t))=\alpha(t)$, $E_t(\beta(t))=\beta(t)$. Let $x\in\a$. Differentiating the identity 
 $$
 <E_t(\alpha(t)),x>=<\alpha(t),x>
 $$
 one obtains
 $$
 <dE_t(\alpha(t)),x>+<E_t(\dot{\alpha}(t)),x>=<\dot{\alpha}(t),x>.
 $$
 This last term equals $<[dE_t,E_t](\alpha(t)),x>$. Note that  
 $$
 E_t(dE_t(\alpha(t)))=E_t(dE_t(E_t(\alpha(t))))=0.
 $$
 Therefore 
 $$
 <[dE_t,E_t](\alpha(t)),x>=<dE_t(\alpha(t)),x>.
 $$
 Then $<E_t(\dot{\alpha}(t)),x>=0$, i.e. $E_t(\dot{\alpha}(t))=0$. Conversely, if  a map $\gamma(t)$ takes values in $\b_t$ and verifies $E_t(\dot{\gamma}(t))=0$, then it is a solution of the transport equation.
 
 The curve $\alpha(t)\beta(t)$ takes values in $\b_t$. 
   Also it is clear that the product rule applies for the derivative of $\alpha(t)\beta(t)$ (as they are $\|\ \|_\infty$ uniformly bounded on closed bounded intervals). Then
 $$ E_t(\frac{d}{dt}(\alpha(t)\beta(t)))=E_t(\dot{\alpha}(t)\beta(t))+E_t(\alpha(t)\dot{\beta}(t))=E_t(\dot{\alpha}(t))\beta(t)+\alpha(t)E_t(\dot{\beta}(t))=0,
 $$
 i.e. $\alpha(t)\beta(t)$ is a solution of the transport equation, with initial condition $ab$. It follows that
 $$
 \theta_t(ab)=G_t(ab)=\alpha(t)\beta(t)=\theta_t(a)\theta_t(b). 
 $$  
 \end{proof}
 
 It was shown above that a solution that starts in $R(E_0)=\b_0$, remains in $R(E_t)=\b_t$ at time $t$. The intertwining identity implies that the same is true for the kernels: if $E_0(a)=0$, then $E_t(\alpha(t))=0$.
In other words, if $a\in\a$ is decomposed as
$$
a=b+z \, \ \ b\in\b_0 \hbox{ and } E_0(z)=0,
$$ 
putting $\beta(t)=G_t(b)$ and $z(t)=G_t(z)$ the solutions with initial conditions $b$ and $z$, then
$$
\alpha(t)=\beta(t)+z(t) \, \ \ \beta(t)\in\b_t \hbox{ and } E_t(z(t))=0,
$$
which is an orthogonal decomposition. The next result shows that their derivatives are also orthogonal for all $t$, though the role of the subspaces is reversed.
\begin{prop}
With the above notations, $E_t(\dot{\beta}(t))=0$ and $\dot{z}(t)\in\b_t$
\end{prop}
\begin{proof}
 As it was shown in the proof of the previous theorem, the solution $\beta(t)$ verifies $\dot{\beta}(t)=dE_t(\beta(t))$, as well as $E_t(dE_t(\beta(t)))=0$. Putting these two together gives $E_t(\dot{\beta}(t))=0$. 
 
On the other hand, since $E_t(z(t))=0$, 
$$
\dot{z}(t)=[dE_t,E_t](z(t))=E_t(dE_t(z(t))),
$$
i.e. $\dot{z}(t)\in \b_t$.
\end{proof}

 \section{Systems of projections}
 
 Let $\p=(p_1,p_2,\dots)$ be a (finite or infinite) system of projections in $\a$, i.e. a sequence of pairwise orthogonal projections which strongly sum $1$. Such a system gives rise to a conditional expectation: 
 $$
 E_\p:\a\to \b\subset \a , \ \ E_\p(x)=\sum_{i\ge 1} p_ixp_i.
 $$
 The range of this conditional expectation is the sub-algebra $\b$ of elements of $\a$ which commute with all $p_i$, $i\ge 1$. Suppose that a curve $\p(t)=(p_1(t),p_2(t),\dots)$, $t\in I$ of systems of projections is given, and that it satisfies that
 $$
 I\ni t\mapsto p_i(t)\xi \in\h
 $$
 is $C^1$ for all $\xi\in\h$ and every $i\ge 1$. We shall examine  the meaning of the  smoothness condition on the curve $E_t=E_{\p(t)}$. We show that if $t\mapsto E_t(a)\xi$ is continuously differentiable (for any $a\in\a$ and $\xi\in\h$), then  Hypothesis (\ref{hipotesis}) holds.
 
Our first elementary observation is that if the system is finite, then these conditions are  fulfilled.

\begin{prop}
Suppose that the system $\p(t)$ is finite, i.e. $\p(t)=(p_1(t),\dots,p_n(t))$, and that for each $j=1,\dots, n$, the curve $p_j(t)\xi$ is $C^1$ in $\h$. Then curve $E_t$ verifies that
$E_t(a)\xi$ is $C^1$ in $\h$ for each $a\in\a$ and $\xi\in\h$, and $dE_t$ is bounded in $\h$.
\end{prop}
\begin{proof}
Pick $a\in\a$ and $\xi\in\h$. Then $E_t(a)\xi$ is $C^1$. Indeed, a straightforward computation shows that the product rule holds and that 
$$
\frac{d}{d t} E_t(a)\xi=\sum_{i=1}^n \dot{p}_i(t)ap_i(t)\xi+p_i(t)a\dot{p}_i(t)\xi.
$$
This map is clearly continuous. Next note that for each $j$, the map $\xi\mapsto \dot{p}_j(t)\xi$ is linear and everywhere defined in $\h$. Moreover, it is symmetric:
$$
<\dot{p}_j\xi,\eta>=\frac{d}{d t}<p_j(t)\xi,\eta>=\frac{d}{d t}<\xi,p_j(t)\eta>=<\xi,\dot{p}_j(t)\eta>.
$$
Therefore, by the closed graph theorem, it is a bounded operator. Since it is defined as a strong limit, it takes values in $\a$, i.e. $\dot{p}_j\in\a$. The operator $dE_t$  coincides in $\a$ with
$$
\sum_{i=1}^n L_{\dot{p}_i(t)}R_{p_i(t)}+L_{p_i(t)}R_{\dot{p}_i(t)},
$$
which is clearly bounded (Here $L_a,R_a$ denote left and right multiplication by $a\in\a$). Moreover, by the uniform boundedness principle, for $t\in J\subset I$, a closed a bounded sub-interval, the norms $\|\dot{p}_j(t)\|_\infty$ are uniformly bounded by $C$ (which can be chosen independent of $j$ as well). Therefore it is apparent that $dE_t$ is bounded in $\h$: 
$$
\|dE_t(a)\|_2\le nC\|a\|_2, \  t\in J.
$$
\end{proof}
We  restrict now to  infinite systems.
First we discuss a condition which implies the regularity of the curve $E_t$. Namely the following, which was studied in \cite{weaksmooth} for expectations in the algebra of compact operators.

\begin{defi}
We shall say that the curve of systems of projections $\p(t)$ has square summable derivatives if for every closed bounded subinterval $J\subset I$, there exists a constant $D_J$ such that
\begin{equation}\label{hipotesisregularidad}
\sum_{i\ge 1} \|\dot{p}_i(t)\xi\|_2^2\le D_J \|\xi\|_2^2
\end{equation} 
for every $\xi\in\h$ and $t\in J$. 
\end{defi}
\begin{prop}
The curve $\p(t)$ has square summable derivatives  (\ref{hipotesisregularidad})  if and only if there exists a strongly $C^1$ curve $u_t$, $t\in I$,  of unitary operators in $\a$ such that $p_i(t)=u_tp_i(0)u_t^*$ for all $i\ge 1$.
\end{prop}
\begin{proof}
Suppose first that inequality (\ref{hipotesisregularidad}) holds. Then we claim that for any $\xi\in\h$ the series 
$$
\sum_{i\ge 1}p_i(t)\dot{p}_i(t)\xi
$$
is convergent in $\h$. Indeed, note that since the vectors $p_i(t)\dot{p}_i(t)\xi$ are pairwise orthogonal,
$$
\|\sum_{i\ge N+1}p_i(t)\dot{p}_i(t)\xi\|_2^2= \sum_{i\ge N+1}\|p_i(t)\dot{p}_i(t)\xi\|_2^2\le \sum_{i\ge N+1}\|\dot{p}_i(t)\xi\|_2^2,
$$
which tends to $0$ as $N$ goes to $\infty$. Then this series produces an everywhere defined linear operator
$$
\Delta_t\xi=\sum_{i\ge 1}p_i(t)\dot{p}_i(t)\xi.
$$
This operator has an everywhere defined adjoint, given by the series
$$
\Delta^*\xi=\sum_{i\ge 1}\dot{p}_i(t)p_i(t)\xi,
$$
which is weakly convergent in $\h$:
$$
<\Delta^*\xi,\eta>=\sum_{i\ge 1}<\dot{p}_i(t)p_i(t)\xi,\eta>=\sum_{i\ge 1}<\xi,p_i(t)\dot{p}_i(t)\eta>.
$$
Therefore, by the closed graph theorem, $\Delta_t$ is bounded, and since it is defined as a strong limit of elements of $\a$, $\Delta_t\in\a$. Note that the identity $\dot{p}_i(t)=\dot{p}_i(t)p_i(t)+p_i(t)\dot{p}_i(t)$ implies
that, since $\sum_{i\ge 1}p_i(t)\xi=\xi$ and this series converges uniformly in closed bounded sub-intervals,
$$
0=\frac{d}{d t}\sum_{i\ge 1}<p_i(t)\xi,\eta>=\sum_{i\ge 1}<\dot{p}_i(t)\xi,\eta>=\sum_{i\ge 1}<\dot{p}_i(t)p_i(t)+p_i(t)\dot{p}_i(t)\xi,\eta>
$$
$$
=<\Delta_t^*\xi+\Delta_t\xi,\eta>,
$$
i.e. $\Delta_t$ is anti-hermitic. Furthermore, the hypothesis that the curve $\p(t)$ has square summable derivatives (\ref{hipotesisregularidad}), implies that on closed bounded sub-intervals, the series that defines $\Delta_t$ is uniformly convergent. Therefore the map
$$
I\ni t\mapsto \Delta_t\xi\in \h
$$
is continuous, that is $t\mapsto \Delta_t\in\a$ is strongly continuous.
For any $\xi_0\in \h$, consider the linear differential equation in $\h$
\begin{equation}\label{ecuaciondelta}
\left\{ \begin{array}{l} \dot{\mu}(t)=-\Delta_t\mu(t) \\ \mu(0)=\xi_0. \end{array}\right.
\end{equation}
It was shown in \cite{weaksmooth} in a different context, that the unitary propagator $u_t$ of this equation, (defined by $u_t\xi_0=\mu(t)$),  verifies 
$$
u_tp_i(0)u_t^*=p_i(t), \ \ i\ge 1.
$$
The computation is formally identical in this context, and thus these relations hold. Moreover, apparently $u_t\in\a$, and the map $t\mapsto u_t\xi_0$ is $C^1$ for every $\xi_0\in \h$.

Conversely, suppose the existence of a strongly $C^1$ curve  $u_t$ of unitaries in $\a$ such that $u_tp_i(0)u_t^*=p_i(t)$ for $i\ge 1$. Then the product rule holds and
$$
\dot{p}_i(t)\xi=\dot{u}_t p_i(0)u_t^*\xi+u_tp_i(0)\dot{u}_t^*\xi.
$$
Then $\|\dot{p}_i(t)\xi\|_2\le \|\dot{u}_t p_i(0)u_t^*\xi\|_2+\|p_i(0)\dot{u}_t^*\xi\|_2$.
Note that for any closed bounded subinterval $J\subset I$,  the family of vectors $\{\dot{u}_t\xi: t\in J\}$ is uniformly bounded. Therefore, by the uniform boundedness principle, $\|\dot{u}_t\|\le K_J$ for all $t\in J$. Then, using that $p_i(0)$ are pairwise orthogonal and sum $1$,
$$
\sum_{i\ge 1} \|\dot{u}_t p_i(0)u_t^*\xi\|_2^2\le K_J^2\sum_{i\ge 1} \|p_i(0)u_t^*\xi\|_2^2=K_J^2\|u_t^*\xi\|_2^2=K_J^2\|\xi\|_2^2,
$$
and
$$
\sum_{i\ge 1}\|p_i(0)\dot{u}_t^*\xi\|_2^2=\|\dot{u}_t^*\xi\|_2^2\le K_J^2\|\xi\|_2^2.
$$
Then 
$$
\sum_{i\ge 1} \|\dot{p}_i(t)\xi\|_2^2\le 4K_J^2\|\xi\|_2^2,
$$
for $t\in J$.
\end{proof}
\begin{rem}
Note that (under the assumption (\ref{hipotesisregularidad}) that the system of projections has square summable derivatives), the unitaries $u_t$ provide anoter way to intertwine $E_0$ and $E_t$. Indeed, put $\Omega_t=Ad(u_t)$ ($\Omega_t(x)=u_txu_t^*$), then
$$
\Omega_tE_0\Omega_t^{-1}(x)=u_t\sum_{i\ge 1} u_tp_i(0)u_t^*xu_tp_i(0)u_t^*=\sum_{i\ge 1} p_i(t)xp_i(t)=E_t(x).
$$
We shall consider the relation between $\Omega_t$ and $G_t$ below. Our purpose now is to use this inner automorphisms to prove the regularity of the curve $E_t$. To this effect, note that for each $a\in\a$ and $\xi\in\h$, the map 
$I\ni t\mapsto \Omega_t(a)\xi$ is $C^1$. Indeed,
$$
\frac1h\{u_{t+h}au^*_{t+h}\xi-u_tau^*_t\xi\}=\frac1h\{u_{t+h}a(u^*_{t+h}\xi-u^*_t\xi)\}+\frac1h\{u_{t+h}au^*_t\xi-u_tau^*_t\xi\}.
$$
The second term tends to $\dot{u}_tau^*_t\xi$ as $h\to 0$, because $u_t$ is strongly $C^1$. The first term tends to $u_ta\dot{u}_t^*\xi$. Indeed, $\|\frac1h\{u_{t+h}a(u^*_{t+h}\xi-u^*_t\xi)\}-u_ta\dot{u}_t^*\xi\|_2$ is bounded by
$$
 \|u_{t+h}a\frac1h\{u^*_{t+h}\xi-u^*_t\xi\}-u_{t+h}a\dot{u}_t^*\xi\|_2+\|u_{t+h}a\dot{u}_t^*\xi-u_ta\dot{u}_t^*\xi\|_2
$$
$$
\le
\|a\frac1h\{u^*_{t+h}\xi-u^*_t\xi\}-a\dot{u}_t^*\xi\|_2+\|u_{t+h}\eta-u_t\eta\|_2,
$$
where $\eta=a\dot{u}_t^*\xi$. Clearly both terms tend to $0$. Finally, the derivative of $\Omega_t(a)\xi$ equals
$$
\dot{\Omega}_t(a)\xi=\dot{u}_tau_t^*\xi+u_ta\dot{u}_t^*\xi,
$$
which is clearly continuous.
 \end{rem}
 Next we show that condition (\ref{hipotesisregularidad}) guarantees  that equation (\ref{trasporte}) has existence and uniqueness of solutions.
\begin{prop}
If the system of projections $\p(t)$ has square summable derivatives (\ref{hipotesisregularidad}), then the map $I\ni t\mapsto E_t(a)\xi\in \h$ is $C^1$. Moreover, the derivative $dE_t$ extends to a bounded operator in $\h$.
\end{prop}

\begin{proof}
As seen above, $E_t(x)=\Omega_t(E_0(\Omega_t^{-1}(x)))$. Note that for each $x\in\a$, both $\Omega_t(x)$ and $\Omega_t^{-1}(x)=u_t^*xu_t$ are strongly $C^1$. Then for each $x\in\a$ and $\xi\in\h$,
$$
\frac1h\{E_{t+h}(x)\xi-E_t(x)\xi\}=\Omega_{t+h}E_0(\frac1h\{\Omega_{t+h}(x)-\Omega_t(x)\})\xi+\frac1h\{(\Omega_{t+h}(x)\eta-\Omega_t(x)\eta\},
 $$
where $\eta=E_0(\Omega_t^{-1}(x))\xi$. The first term tends to $\Omega_tE_0\dot{\Omega}_t(x)\xi$: put 
$$
b_h=E_0(\frac1h\{\Omega_{t+h}(x)-\Omega_t(x)\}),
$$
which tends strongly to $b_0=E_0(\dot{\Omega}_t(x))$ (because $E_0$ is strongly continuous), then
$$
\|\Omega_{t+h}(b_h)\xi-\Omega_t(b_0)\xi\|_2\le \|\Omega_{t+h}(b_h)\xi-\Omega_t(b_h)\xi\|_2+\|\Omega_{t}(b_h)\xi-\Omega_t(b_0)\xi\|_2.
$$
The second term clearly tends to $0$. The first term is bounded by 
$$
\|u_{t+h}b_h(u^*_{t+h}-u^*_t)\xi\|_2+\|u_tb_h(u^*_{t+h}-u_t^*)\xi\|_2\le 2\|b_h\|_\infty\|u^*_{t+h}\xi-u_t^*\xi\|_2.
$$
This term tends to zero because   the involution $*$ is strongly continuous ($\a$ is finite) and $\|b_h\|_\infty$ is bounded for $|h|$ small:
$$
\|b_h\|_\infty \le \|\frac1h\{\Omega_{t+h}(x)-\Omega_t(x)\}\|_\infty,
$$
with $\frac1h\{\Omega_{t+h}(x)-\Omega_t(x)\}$ strongly convergent, and therefore locally $\| \ \|_\infty$-bounded.

Note that, in the above notations, $\xi\mapsto \dot{u}_t\xi$ is an everywhere defined operator. Clearly $u_t^*\dot{u}_t$ is anti-hermitian:
$$
0=\frac{d}{dt}<u_t\xi,u_t\eta>=<u_t^*\dot{u}_t\xi,\eta>+<\xi,u_t^*\dot{u}_t\eta>.
$$
Then, by the closed graph theorem, $u_t^*\dot{u}_t$ is bounded, and therefore $\dot{u}_t$ is bounded. Also it is clear that, being a strong limit of operators in $\a$, it belongs to $\a$. Then
$$
\dot{\Omega}_t=L_{\dot{u}_t}R_{u_t^*}+R_{\dot{u}_t}L_{u_t^*}
$$
is bounded. Also it is clear that $\Omega_t^{-1}=Ad(u_t^*)$ has the same properties. Then 
$$
dE_t=\dot{\Omega}_tE_0\Omega_t^{-1}+\Omega_tE_0\dot{\Omega^{-1}}_t
$$
is bounded in $\h$.
 \end{proof}
\begin{rem}
In \cite{weaksmooth}, similar results were obtained for the algebra $\k(\h)$ of compact operators. For instance it was shown that if the systems $\p(t)$ consist of more than two projectors, then $\Omega_t$ and $G_t$ differ. It was also shown that they coincide if the system consists of two projections, and that $\Omega_t$ and $G_t$ coincide in $\b_0$. In other words, always under the assumption that inequality (\ref{hipotesisregularidad}) holds, the unitaries $u_t$ of $\a$ which solve equation (\ref{ecuaciondelta}), implement the automorphism $\theta_t$:
$$
\theta_t=Ad(u_t)|_{\b_0}:\b_0\to \b_t.
$$
 We refer the reader to \cite{weaksmooth} for the proofs of these facts, which though  perfomed in $\k(\h)$, are formally identical in our situation.
\end{rem}

We now show that for this class of conditional expectations, given by a system of projections, smoothness of the curve $E_t$ implies Hypothesis (\ref{hipotesis}).

\begin{prop}
Let $\p(t)$, $t\in I$, be a system of projectors and $E_t$ as above, verifying that  $I\ni t \mapsto E_t(a)\xi\in \h$ is $C^1$, for every $a\in\a$ and $\xi\in\h$, and that for each $j\ge 1$, $t\mapsto p_j(t)\xi$ is $C^1$.
Then  Hypothesis (\ref{hipotesis}) holds: for each closed and bounded sub-interval $J\subset I$, there exists $C_J$ such that
$$
\int_J\|dE_t(a)\|_2^2 dt\le C_J\|a\|_2^2
$$
for each $a\in\a$.
\end{prop}
\begin{proof}
Note that the map $t\mapsto p_j(t)\in\a$ is a solution of equation (\ref{trasporte}). Since $p_j(t)\in\b_t$, this equation becomes simpler, as seen in the previous section. Namely, one has to show that 
$$
\dot{p}_j(t)=dE_t(p_j(t)).
$$
Indeed:
$$
dE_t(p_j(t))=\sum_{i\ge 1} \dot{p}_i(t)p_j(t)p_i(t)+p_i(t)p_j(t)\dot{p}_i(t)=\dot{p}_j(t)p_j(t)+p_j(t)\dot{p}_j(t)=\dot{p}_j(t),
$$
where the last identity follows from differentiating $p_j(t)p_j(t)=p_j(t)$. 
Then we can bound the operator norm of $\dot{p}_j(t)\in\a$:
$$
\|\dot{p}_j(t)\|_\infty=\|dE_t(p_j(t))\|_\infty\le \|dE_t\|_{\infty,\infty}\le D_J
$$
for a constant $D_J$ independent of $t\in J$. Then
$$
\|\sum_{i\ge 1}\dot{p}_i(t)ap_i(t)\|_2^2=\sum_{i\ge 1}\tau(p_i(t)a^*(\dot{p}_i(t))^2ap_i(t))\le D_J^2 \sum_{i\ge 1} \tau(p_i(t)a^*ap_i(t))
$$
$$
=D_J^2\sum_{i\ge 1} \tau(p_i(t)a^*a)=D_J^2\tau(a^*a)=D_J^2\|a\|_2^2.
$$
Analogously, $\|\sum_{i\ge 1}p_i(t)a\dot{p}_i(t)\|_2^2\le D_J^2\|a\|_2^2$.
Then 
$$
\|dE_t(a)\|_2^2=\|\sum_{i\ge 1}\dot{p}_i(t)ap_i(t)+\sum_{i\ge 1}p_i(t)a\dot{p}_i(t)\|_2^2\le 4D_j^2 \|a\|_2^2.
$$
Therefore
$$
\int_J \|dE_t(a)\|_2^2 d t \le 4|J|D_J^2\|a\|_2^2.
$$
\end{proof}

\bigskip

\noindent

Esteban Andruchow and Gabriel Larotonda\\

Instituto de Ciencias \\
Universidad Nacional de Gral. Sarmiento \\
J. M. Gutierrez 1150 \\
(1613) Los Polvorines \\
Argentina  \\

and \\

Instituto Argentino de Matem\'atica\\
''Alberto P. Calder\'on'', CONICET \\
Saavedra 15, 3er. piso \\
(1083) Buenos Aires \\
Argentina.

e-mails: eandruch@ungs.edu.ar, glaroton@ungs.edu.ar

\end{document}